\pdfoutput=1
\documentclass[colorlinks = true, a4paper, oneside, reqno, 12pt]{amsart}

\usepackage[top = 3cm, bottom = 3cm, left = 2.5cm, right = 2.5cm, marginpar=2cm]{geometry}
\usepackage[T1]{fontenc}
\usepackage{lmodern}

\usepackage{microtype}

\usepackage{booktabs}

\usepackage[style=alphabetic, sorting=nyt, maxbibnames=99]{biblatex}
\usepackage[shortlabels]{enumitem}
\usepackage{calc}

\setlist{topsep=-0.4ex, itemsep=0.2ex, parsep=0ex}

\newlength{\bulletindent}
\newlength{\bulletbodyindent}
\setlist[itemize,1]{
    left=\bulletindent .. \bulletbodyindent
}

\usepackage{mathtools, amssymb}
\usepackage[mathcal]{eucal}

\usepackage{tikz}
\usepackage{tikz-cd}

\usepackage{amsthm}

\newtheorem{theorem}{Theorem}
\newtheorem{proposition}[theorem]{Proposition}

\theoremstyle{definition} 
\newtheorem{definition}[theorem]{Definition}

\theoremstyle{remark}
\newtheorem*{remark}{Remark}

\usepackage{makecell}

\renewcommand{\leq}{\leqslant}
\renewcommand{\geq}{\geqslant}

\usepackage{hyperref}
\hypersetup{allcolors=[rgb]{0.1,0.1,0.4}}

\newcommand{\A}{\mathcal{A}}
\newcommand{\B}{\mathcal{B}}
\newcommand{\C}{\mathcal{C}}
\newcommand{\D}{\mathcal{D}}

\renewcommand{\Phi}{\varPhi}
\renewcommand{\Psi}{\varPsi}
\renewcommand{\Upsilon}{\varUpsilon}
\renewcommand{\Omega}{\varOmega}

\newcommand{\V}{\mathcal{V}}
\newcommand{\W}{\mathcal{W}}
\newcommand{\Set}{\mathrm{\mathcal{S}et}}
\newcommand{\Cat}{\mathrm{\mathcal{C}at}}
\newcommand{\Lens}{\mathrm{\mathcal{L}ens}}

\newcommand{\Fam}{\mathrm{\mathcal{F}am}}
\newcommand{\wSet}{\mathrm{w\mathcal{S}et}}

\renewcommand{\Vec}{\mathrm{\mathcal{V}ec}}

\newcommand{\COF}{\mathrm{\mathbb{C}of}}

\mathchardef\mhyphen="2D
\newcommand{\VCat}{\V\mhyphen\Cat}
\newcommand{\VCOF}{\V\mhyphen\COF}
\newcommand{\VLens}{\V\mhyphen\Lens}

\DeclareMathOperator{\cod}{cod}
\DeclareMathOperator{\id}{id}

\DeclareMathOperator{\obj}{obj}

\DeclareMathOperator{\flt}{flat}
\DeclareMathOperator{\dist}{d}
\DeclareMathOperator{\disc}{disc}

\makeatletter
\newcommand{\pto}{}% just for safety
\newcommand{\pgets}{}% just for safety

\DeclareRobustCommand{\pto}{\mathrel{\mathpalette\p@to@gets\to}}
\DeclareRobustCommand{\pgets}{\mathrel{\mathpalette\p@to@gets\gets}}

\newcommand{\p@to@gets}[2]{%
  \ooalign{\hidewidth$\m@th#1\mapstochar\mkern5mu$\hidewidth\cr$\m@th#1\to$\cr}%
}
\makeatother

\usetikzlibrary{decorations.markings}
\usetikzlibrary{calc}
\usetikzlibrary{fit}

\newcommand{\arrowLength}{2.3pt}
\newcommand{\arrowWidth}{4.6pt}

\tikzcdset{
    postmark/.style={
        postaction={
            decorate,
            decoration={
                markings,
                mark={#1}
            }
        }
    },
    proarrow/.style={
        postmark=at position 0.5 with {
            \draw[-] (0,2pt) -- (0,-2pt);
        },
        labels={
            inner sep=0.75ex,
        }
    },
	morphismTip/.tip={
		Triangle[length=\arrowLength, width=\arrowWidth]
	},
	lensTip/.tip={
		Diamond[length=2*\arrowLength, width=\arrowWidth]
	},
	mark/.style={
		postaction={
			decorate,
			decoration={
				markings,
				mark={#1}
			}
		}
	},
	wideMark/.style={
		labels={
			inner sep=0.8ex,
		}
	},
	functor/.style={
		wideMark,
		mark=at position {0.5*\pgfdecoratedpathlength +
		\arrowLength} with {\arrow{morphismTip}},
	},
	cofunctor/.style={
		wideMark,
		mark=at position {0.5*\pgfdecoratedpathlength -
		\arrowLength} with {\arrowreversed{morphismTip}},
	},
	lens/.style={
		wideMark,
		mark=at position {0.5*\pgfdecoratedpathlength +
		\arrowLength} with {\arrow{lensTip}},
	},
	inline/.style={
		cramped,
		sep=small,
	},
	display style/.style={
        cells={font=\everymath\expandafter{\the\everymath\displaystyle}},
    },
	split arrows style/.style={
		display style,
		cramped,
		row sep=tiny,
	},
	arrows={
		-{Computer Modern Rightarrow[length=3.5pt, width=7pt]},	
		semithick
	},
    cong/.style={
        "\cong",
        sloped,
        phantom
    },
    short equals/.style={
        "=",
        sloped,
        phantom
    },
    isomorphism/.style={
        "\sim"{sloped, inner sep=0.08ex},
        every label/.append style={inner sep=1.4ex}
    },
    isomorphism swap/.style={
        "\sim"{sloped, swap, inner sep=0.25ex},
        every label/.append style={inner sep=1.4ex}
    },
}

\begin{document}

\title{An introduction to enriched cofunctors}

\author{Bryce Clarke}

\address{Inria Saclay Centre \\ Palaiseau, France}

\email{bryce.clarke@inria.fr}

\author{Matthew Di Meglio}

\address{University of Edinburgh \\ Edinburgh, Scotland}

\email{m.dimeglio@ed.ac.uk}



%----------------------------------------------------------------------%
% Abstract                                                             %
%----------------------------------------------------------------------%
\begin{abstract}
Cofunctors are a kind of map between categories which lift 
morphisms along an object assignment. 
In this paper, we introduce cofunctors between categories enriched in a 
distributive monoidal category. 
We define a double category of enriched categories, enriched functors, 
and enriched cofunctors, whose horizontal and vertical $2$-categories 
have $2$-cells given by enriched natural transformations between 
functors and cofunctors, respectively. 
Enriched lenses are defined as a compatible enriched functor and 
enriched cofunctor pair; 
weighted lenses, which were introduced by Perrone, are precisely lenses 
enriched in weighted sets. 
Several other examples are also studied in detail.
\end{abstract}

\maketitle

%----------------------------------------------------------------------%
% Section 1: Introduction                                              %
%----------------------------------------------------------------------%
\section{Introduction}
\label{sec:introduction}

Cofunctors were introduced by Aguiar \cite{Agu97} as a kind of morphism 
between categories, generalising the notion of \emph{comorphism} 
between Lie groupoids \cite{HM93}. 
Although cofunctors have not played a significant role in classical 
category theory, there has recently been a growing interest in their 
properties and applications.
Ahman and Uustalu proved that polynomial comonads on $\Set$ are 
precisely small categories \cite{AU16}, whilst the comonad morphisms 
correspond to cofunctors \cite{AU17}. 
The category $\Cat^{\sharp}$ of small categories and cofunctors, and 
its relationship to polynomial functors, is the subject of ongoing 
work of Spivak and Nui \cite{ SN22}.
Applications of cofunctors in mathematics and computer science have
also appeared in a number of places 
\cite{Gar19, Cla20b, Spi20, CG21, Gar21, Spi22}. 

The purpose of this paper is to generalise the notion of cofunctor to 
the setting of enriched category theory. 
One reason to do this is to provide an enriched counterpart to the 
theory of \emph{internal cofunctors} developed by Aguiar 
\cite[Section~4.2]{Agu97}.
The goal of both enriched and internal category theory is to not 
only reach a better understanding of (ordinary) category theory, 
but to also extend known concepts to a wider range of examples. 
The primary contributions of this paper are to establish the basic 
theory of cofunctors enriched over a 
\emph{distributive monoidal category}, and to show how this captures 
a number of examples in the literature 
(see Section~\ref{sec:examples}).

A key motivation for initiating the study of enriched cofunctors 
comes from recent work by Perrone on a variant of delta lenses 
called \emph{weighted lenses}~\cite{Per21}. 
\emph{Delta lenses} are a kind of morphism between categories~\cite{DXC11},
generalising the notion of split opfibration~\cite{JR13}.
A \emph{weighted category} is a category in which every morphism is 
equipped with an extended non-negative real number called a \emph{weight}~\cite{Gra07},
and a weighted lens is the same as 
a delta lens whose underlying functor and cofunctor are 
\emph{weight non-increasing}. 
However, weighted categories are precisely categories enriched in 
weighted sets, so it is natural to wonder whether weighted lenses are examples of some notion of \emph{enriched lens}. We answer this question affirmatively (see Definition~\ref{definition:enriched-lens} and Section~\ref{sec:weighted-lenses}).

Similar to Ahman and Usutalu's characterisation of delta lenses as consisting of a functor 
and cofunctor that are compatible with each other in a certain way~\cite{AU17}, we define enriched lens in terms of an analogous notion of \textit{compatibility} for enriched functors and cofunctors. To fully understand compatibility, it is important, first, to appreciate the duality between the concepts of enriched functor and cofunctor.

At first glance, functors and cofunctors appear to be quite different. 
A functor assigns objects and morphisms in the same 
direction, whereas a cofunctor lifts morphisms in the opposite direction 
to its assignment on objects. 
In addition, enriched functors may be enriched in \emph{any} monoidal category, whilst our definition of enriched cofunctor uses the assumption that the base of enrichment is distributive monoidal in an essential way. There is, however, a hidden duality between functors and cofunctors (see Table~\ref{table:functors-and-cofunctors}). 

For a functor $F \colon \A \rightarrow \B$, the data of the functions
$F_{a, b} \colon \A(a, a') \rightarrow \B(Fa, Fa')$ 
is equivalent, by the universal property of the coproduct, 
to the function~$[F_{a, x}]$ in Table~\ref{table:functors-and-cofunctors}. 
Cofunctors $\Phi \colon \A \pto \B$ 
are dual to functors in the sense that the direction of 
the corresponding morphism assignments are reversed. 
The compatibility between a cofunctor~$\Phi$ and a functor $F$ may then be expressed as
the requirement that their morphism assignments, as depicted in 
Table~\ref{table:functors-and-cofunctors}, 
are a section-retraction pair.
This duality extends to functors and cofunctors enriched in a distributive monoidal category. 

\begin{table}
	\centering
	\setlength\tabcolsep{1em}
	\caption{Comparison of (enriched) functors and cofunctors.}
	\begin{tabular}{ccc}\toprule
		&
		Object assignment&
		Morphism assignment
	\\ \midrule 
		\makecell[c]{functor \\ $F \colon \A \to \B$}
		&
		$F \colon \obj(\A) \rightarrow \obj(\B)$ &
		$[F_{a, x}] \colon \negthickspace\negthickspace \displaystyle\sum_{x \in F^{-1}\{b\}} \A(a, x) 
		\rightarrow \B(Fa, b)$
	\\[3ex]
		\makecell[c]{cofunctor \\ $\Phi \colon \A \pto \B$}&
		$\Phi \colon \obj(\A) \rightarrow \obj(\B)$ &
		$\Phi_{a, b} \colon \B(Fa, b)
		\rightarrow \negthickspace \displaystyle\sum_{x \in \Phi^{-1}\{b\}} \negthickspace 
		\A(a, x)$
	\\ \bottomrule
	\end{tabular}
	\label{table:functors-and-cofunctors}
\end{table}

It is also natural to wonder whether there are more general notions of 
``compatibility'' between (enriched) functors and cofunctors, 
and whether there is a setting in which they may be studied together. One approach is to try to construct a double category of 
(enriched) categories, functors, and cofunctors, an idea
first suggested by Par\'{e} (for the non-enriched setting)~\cite{Par20}.
For (ordinary) functors and cofunctors, an appropriate choice of 
cell for such a double category is the notion of \emph{compatible square}, and delta lenses are compatible squares of a certain form~\cite{Cla22, DiM21}. 
For enriched functors and cofunctors, 
we define an analogous notion of compatible square, yielding the
flat double category $\VCOF_{\flt}$; enriched lenses are likewise certain cells therein. 

It is unsatisfying that we cannot recover from $\VCOF_{\flt}$ the notion of enriched natural transformation. Surprisingly, there is another double category $\VCOF$, which has  
$\VCOF_{\flt}$ as a double subcategory, whose horizontal 
$2$-category is the $2$-category of enriched categories, 
functors, and natural transformations. 
Furthermore, 
the $2$-cells of the vertical $2$-category of $\VCOF$ are a notion of enriched natural 
transformation between enriched cofunctors analogous to Aguiar's \emph{natural cotransformation} in the internal setting. 

%----------------------------------------------------------------------%
% Subsection                                                           %
%----------------------------------------------------------------------%
\subsection*{Outline of the paper}

In Section~\ref{sec:enriched-categories}, 
we recall the definition of the $2$-category $\VCat$ of categories, functors, 
and natural transformations
enriched in a monoidal category $\V$.

In Section~\ref{sec:enriched-cofunctors}, 
we define the new notions of cofunctor enriched in a distributive monoidal category $\V$, and enriched natural transformation between them. Additionally, we construct a $2$-category $\VCat^{\sharp}$
of enriched categories, cofunctors, and natural transformations.

In Section~\ref{sec:double-category}, 
we define the double category $\VCOF$ of enriched categories, functors,
cofunctors, and natural transformations,
and show how $\VCat$ and $\VCat^{\sharp}$ are, respectively, the
horizontal and vertical $2$-categories of this double category.

In Section~\ref{sec:enriched-lenses}, we introduce the notion of enriched lens. We also define a flat double subcategory $\VCOF_{\flt}$ of $\VCOF$
whose cells are compatible squares, 
and characterise enriched lenses as compatible squares of a certain 
form. 

In Section~\ref{sec:examples}, 
we describe the cofunctors and lenses enriched in 
\begin{itemize}
\item the category of sets and functions equipped with the cartesian 
monoidal product;
\item the category of vector spaces equipped with the usual tensor product;
\item the thin category of extended non-negative real numbers equipped 
with addition as the monoidal product;
\item the free coproduct completion of a monoidal category; and,
\item the category of weighted sets. 
\end{itemize}

%----------------------------------------------------------------------%
% Subsection                                                           %
%----------------------------------------------------------------------%
\subsection*{Future outlook} 
The current version of this paper is the first part of a larger body of 
research currently under preparation. 
In a forthcoming version of this paper, 
we plan to
\begin{itemize}
\item explain the equivalence between enriched cofunctors and 
certain spans of enriched functors, which holds when the base of enrichment is extensive
\cite{CDM22};
\item introduce the theory of 
\emph{monad retromorphisms} in a double category with a functorial choice of companions~\cite{Cla22}—in the double categories of matrices in a distributive monoidal category $\V$, of spans in a category $\mathcal{E}$ with pullbacks, and of comodules in a monoidal category $\W$ with nice equalisers, these are, respectively, cofunctors enriched in $\V$, cofunctors internal to $\mathcal{E}$~\cite{Cla20a}, and cofunctors internal to $\W$~\cite{Agu97}; and,
\item give a new perspective on labelled transition systems and bisimulation in terms of enrichment in the free quantale $P(L^*)$ on the free monoid $L^*$ of the label set $L$.
\end{itemize}

%----------------------------------------------------------------------%
% Section 2: Background.                                               %
%----------------------------------------------------------------------%
\section{Enriched categories, functors, and natural transformations}
\label{sec:enriched-categories}

In this section, we recall from Kelly \cite{Kel05} 
the concepts of enriched category, 
enriched functor, and enriched natural transformation over a fixed 
monoidal category $\V$. We also define the corresponding
$2$-category $\VCat$. 

\begin{definition}
\label{definition:enriched-category}
An \emph{enriched category} $\A$ consists of a set $\obj(\A)$ of objects, 
a hom-object $\A(x, y)$ in $\V$ for each pair of objects,
an identity element $\eta_{x} \colon I \rightarrow \A(x, x)$ 
for each object, 
and a composition map 
$\mu_{x, y, z} \colon \A(x, y) \otimes \A(y, z) \rightarrow \A(x, z)$
for each triple of objects, 
such that the following diagrams, 
which encode the unitality and associativity laws, commute. 
\begin{equation*}
\begin{tikzcd}[column sep = huge]
I \otimes \A(x, y) 
\arrow[d, "\eta_{x} \otimes \id"']
\arrow[rd, "\cong"]
& 
& \A(x, y) \otimes I 
\arrow[d, "\id \otimes \eta_{y}"]
\arrow[ld, "\cong"']
\\
\A(x, x) \otimes \A(x, y) 
\arrow[r, "\mu_{x, x, y}"']
& \A(x, y)
& \A(x, y) \otimes \A(y, y) 
\arrow[l, "\mu_{x, y, y}"]
\end{tikzcd}
\end{equation*}
\begin{equation*}
\begin{tikzcd}[column sep = huge]
\big( \A(w, x) \otimes \A(x, y) \big) \otimes \A(y, z)
\arrow[r, "\mu_{w, x, y} \otimes \id"]
\arrow[d, "\cong"']
& \A(w, y) \otimes \A(y, z)
\arrow[dd, "\mu_{w, y, z}"]
\\
\A(w, x) \otimes \big( \A(x, y) \otimes \A(y, z) \big)
\arrow[d, "\id \otimes \mu_{x, y, z}"']
& \\
\A(w, x) \otimes \A(x, z) 
\arrow[r, "\mu_{w, x, z}"']
& \A(w, z)
\end{tikzcd}
\end{equation*}
\end{definition}

\begin{definition}
\label{definition:enriched-functor}
An \emph{enriched functor} $F \colon \A \rightarrow \B$ 
consists of a function 
$F \colon \obj(\A) \rightarrow \obj(\B)$ 
and a family 
$F_{a, a'} \colon \A(a, a') \rightarrow \B(Fa, Fa')$
of morphisms in $\V$, such that the following diagrams, 
which encode preservation of identities and composites, commute.
\begin{equation*}
\begin{tikzcd}[column sep = tiny, row sep = 2.15em]
& I
\arrow[ld, "\eta_{a}"']
\arrow[rd, "\eta_{Fa}"]
& \\
\A(a, a)
\arrow[rr, "F_{a, a}"']
& & \B(Fa, Fa)
\end{tikzcd}
\quad
\begin{tikzcd}[column sep = huge]
\A(a, a') \otimes \A(a', a'')
\arrow[r, "\mu_{a, a', a''}"]
\arrow[d, "F_{a, a'} \otimes F_{a', a''}"']
& \A(a, a'')
\arrow[d, "F_{a, a''}"]
\\
\B(Fa, Fa') \otimes \B(Fa', Fa'')
\arrow[r, "\mu_{Fa, Fa', Fa''}"']
& \B(Fa, Fa'')
\end{tikzcd}
\end{equation*}
\end{definition}

\begin{definition}
\label{definition:enriched-transformation}
An \emph{enriched natural transformation} $\tau \colon F \Rightarrow G \colon \A \to \B$ of enriched functors consists of a family $\tau_{a} \colon I \rightarrow \B(Fa, Ga)$
of morphisms in $\V$, such that the following diagram, which encodes
naturality, commutes. 
\begin{equation*}
\begin{tikzcd}[column sep={4.4em}]
I \otimes \A(a, a')
\arrow[d, "\tau_{a} \otimes G_{a, a'}"']
& \A(a, a')
\arrow[l, "\cong"']
\arrow[r, "\cong"]
& \A(a, a') \otimes I
\arrow[d, "F_{a, a'} \otimes \tau_{a'}"]
\\
\B(Fa, Ga) \otimes \B(Ga, Ga')
\arrow[r, "\mu_{Fa, Ga, Ga'}"']
&
\B(Fa, Ga')
& \B(Fa, Fa') \otimes \B(Fa', Ga')
\arrow[l, "\mu_{Fa, Fa', Ga'}"]
\end{tikzcd}
\end{equation*}
\end{definition}

There is a $2$-category $\VCat$ whose objects are enriched categories,
whose $1$-cells are enriched functors, 
and whose $2$-cells are enriched natural transformations. 

%----------------------------------------------------------------------%
% Section 2: Main definitions                                          %
%----------------------------------------------------------------------%
\section{Enriched cofunctors}
\label{sec:enriched-cofunctors}

In this section, we introduce the notions of enriched cofunctor and 
enriched natural transformation between enriched cofunctors
over a fixed distributive monoidal
category $\V$. We also define the corresponding $2$-category $\VCat^{\sharp}$. 

\begin{definition}
\label{definition:enriched-cofunctor}
An \emph{enriched cofunctor} $\Phi \colon \A \pto \B$ 
consists of a function
\[\Phi \colon \obj(\A) \rightarrow \obj(\B)\]
called the \emph{object map}, 
and a family
\begin{equation*}
	\Phi_{a, b} \colon 
	\B(\Phi a, b) \longrightarrow 
	\sum_{x \in X} \A(a, x)
	\qquad\qquad
	X = \Phi^{-1}\{b\}
\end{equation*}
of morphisms in $\V$, called \emph{lifting maps}, 
such that the following diagrams, which encode
preservation of identities and composites, commute.
\begin{equation*}
\begin{tikzcd}[column sep = large]
I 
\arrow[r, "\eta_{\Phi a}"]
\arrow[d, "\eta_{a}"']
& \B(\Phi a, \Phi a) 
\arrow[d, "\Phi_{a, \Phi a}"]
\\
\A(a, a)
\arrow[r, "\iota_{a}"']
& \displaystyle\sum_{x \in X} \A(a, x)
\end{tikzcd}
\qquad \qquad X = \Phi^{-1}\{\Phi a\} 
\end{equation*}
\begin{equation*}
\begin{tikzcd}[column sep = huge]
\B(\Phi a, b) \otimes \B(b, b')
\arrow[r, "\mu_{\Phi a, b, b'}"]
\arrow[d, "\Phi_{a, b} \otimes \id"']
& \B(\Phi a, b')
\arrow[dddd, "\Phi_{a, b'}"] 
\\
\big( \displaystyle\sum_{x \in X} \A(a, x) \big) \otimes \B(b, b')
\arrow[d, "\cong"']
& \\
\displaystyle\sum_{x \in X} \A(a, x) \otimes \B(\Phi x, b') 
\arrow[d, "\sum \id \otimes \Phi_{x, b'}"']
& \\
\displaystyle\sum_{x \in X} \A(a, x) \otimes 
\big( \displaystyle\sum_{y \in Y} \A(x, y) \big)
\arrow[d, "\cong"']
& \\
\displaystyle\sum_{y \in Y} \displaystyle\sum_{x \in X} 
\A(a, x) \otimes \A(x, y) 
\arrow[r, "{\sum [\mu_{a, x, y}]}"'] 
& \displaystyle\sum_{y \in Y} \A(a, y)
\end{tikzcd}
\qquad \qquad
\begin{aligned}
X &= \Phi^{-1}\{b\} \\
Y &= \Phi^{-1}\{b'\}
\end{aligned}
\end{equation*}
\end{definition}

There are two surprises in the above definition---the distributivity requirement on the base of enrichment $\V$, 
and the indexing of the sum in the codomain of $\Phi_{a, b}$ 
by the fibre $\Phi^{-1}\{b\}$ rather than all of $\obj(\A)$. 
Both of these are essential for formulating the axiom expressing the compatibility of the cofunctor with composition.

The composite
of enriched cofunctors $\Phi \colon \A \pto \B$ and $\Psi \colon \B \pto \C$ has object map given by the composite of their object maps and lifting maps given by the composite morphism
\begin{gather*}
\C(\Psi\Phi a, c)
	\xrightarrow{\,\Psi_{\Phi a, c}\;}
\sum_{x \in X} \B(\Phi a, x)
	\xrightarrow{\,\sum \Phi_{a, x}\;}
\sum_{x \in X}\sum_{y \in Y_{x}} \A(a, y)
	\xrightarrow{\,\cong\;}
\sum_{z \in Z} \A(a, z).
\\[1em]
X = \Psi^{-1}\{c\} \qquad \qquad 
Y_{x} = \Phi^{-1}\{x\} \qquad\qquad 
Z = (\Psi\Phi)^{-1}\{c\}
\end{gather*}
Whilst it is clear that composition of enriched cofunctors
is strictly unital and also associative up to isomorphism; the proof that associativity is actually \emph{strict} is more subtle.

\begin{definition}
\label{definition:enriched-cotransformation}
An \textit{enriched natural transformation} $\tau \colon \Phi \Rightarrow \Psi \colon \A \pto \B$ of enriched cofunctors consists of a family
\begin{equation*}
\tau_{a} \colon I \longrightarrow 
\sum_{x \in X}\A(a, x)\qquad\qquad X=\Psi^{-1}\{\Phi a\}
\end{equation*}
of morphisms in $\V$, called the \emph{components} of $\tau$, such that the following diagram, which encodes naturality, commutes.
\begin{gather*}
\begin{tikzcd}[column sep=large, ampersand replacement=\&]
I \otimes \B(\Phi a, b)
\arrow[d, "\tau_{a} \otimes \id"']
\& \B(\Phi a, b)
\arrow[l, "\cong"']
\arrow[r, "\cong"]
\& \B(\Phi a, b) \otimes I
\arrow[d, "\Phi_{a, b} \otimes \id"]
\\
\big( \displaystyle\sum_{x \in X} \A(a, x) \big) \otimes \B(\Phi a, b)
\arrow[d, "\cong"']
\& \& \big( \displaystyle\sum_{y \in Y} \A(a, y) \big) \otimes I
\arrow[d, "\cong"]
\\
\displaystyle\sum_{x \in X} \A(a, x) \otimes \B(\Psi x, b)
\arrow[d, "\sum \id \otimes \Psi_{x, b}"']
\& \& \displaystyle\sum_{y \in Y} \A(a, y) \otimes I
\arrow[d, "\sum \id \otimes \tau_{y}"]
\\
\displaystyle\sum_{x \in X} \A(a, x) \otimes 
\big( \displaystyle\sum_{z \in Z} \A(x, z) \big)
\arrow[d, "\cong"']
\& \& \displaystyle\sum_{y \in Y} \A(a, y) \otimes 
\big( \displaystyle\sum_{z \in Z} \A(y, z) \big)
\arrow[d, "\cong"]
\\
\displaystyle\sum_{z \in Z} 
\displaystyle\sum_{x \in X} \A(a, x) \otimes \A(x, z)
\arrow[r, "{\sum [\mu_{a, x, z}]}"']
\&
\displaystyle\sum_{z \in Z} \A(a, z) 
\& \displaystyle\sum_{z \in Z} 
	\displaystyle\sum_{y \in Y} \A(a, y) \otimes \A(y, z)
\arrow[l, "{\sum [\mu_{a, y, z}]}"]
\end{tikzcd}
\\[1em]
X = \Psi^{-1}\{\Phi a\} \qquad \qquad 
Y = \Phi^{-1}\{b\} \qquad \qquad
Z = \Psi^{-1}\{b\} 
\end{gather*}
\end{definition}

\begin{remark}
Aguiar uses the term \emph{natural cotransformation} \cite{Agu97} for 
what we simply refer to as \emph{(enriched) natural transformation}. It will always be clear from context whether a given natural transformation is between a pair of functors or a pair of cofunctors.
\end{remark}

There is a $2$-category $\VCat^{\sharp}$ whose 
objects are enriched categories,
whose $1$-cells are enriched cofunctors, 
and whose $2$-cells are enriched natural transformations. To define the composition of $2$-cells in this $2$-category, it suffices to define the vertical composition of $2$-cells and also the left and right whiskering of $2$-cells with $1$-cells.

The \textit{left whiskering}
\[
\begin{tikzcd}
	\D
		\arrow[r, "\Upsilon", proarrow]&
	\A
		\arrow[r, "\Phi"{name=A}, bend left=40, proarrow]
		\arrow[r, "\Psi"{swap, name=B}, bend right=40, proarrow]
		\arrow[from=A,to=B, phantom, "\Downarrow", "\scriptstyle\tau"{anchor=west, inner sep=0.8ex}]&
	\B
\end{tikzcd}
\]
in $\VCat^{\sharp}$ has component at $d \in \obj(\D)$ given by the composite morphism
\begin{gather*}
I
	\xrightarrow{\,\tau_{\Upsilon d}\;}
\sum_{x \in X} \A(\Upsilon d, x)
	\xrightarrow{\,\sum \Upsilon_{d, x}\;}
\sum_{x \in X} \sum_{y \in Y_{x}} \D(d, y)
	\xrightarrow{\,\cong\;}
\sum_{z \in Z} \D(d, z).
\\[1em]
X = \Psi^{-1}\{ \Phi\Upsilon d\} \qquad \qquad
Y_{x} = \Upsilon^{-1}\{x\} \qquad \qquad
Z = (\Psi\Upsilon)^{-1}\{\Phi\Upsilon d\}
\end{gather*}

The \textit{right whiskering} 
\[
\begin{tikzcd}
	\A
		\arrow[r, "\Phi"{name=A}, bend left=40, proarrow]
		\arrow[r, "\Psi"{swap, name=B}, bend right=40, proarrow]
		\arrow[from=A,to=B, phantom, "\Downarrow", "\scriptstyle\tau"{anchor=west, inner sep=0.8ex}]&
	\B	\arrow[r, "\Omega", proarrow]&
	\C
\end{tikzcd}
\]
in $\VCat^{\sharp}$ has component at $a \in \obj(\A)$ given by the composite morphism
\begin{gather*}
I
	\xrightarrow{\,\tau_{a}\;}
\sum_{x \in X} \A(a, x)
	\xrightarrow{\,[\iota_{x}]\;}
\sum_{y \in Y} \A(a, y).
\\[1em]
X = \Psi^{-1}\{\Phi a\} \qquad \qquad 
Y = (\Omega \Psi)^{-1}\{ \Omega \Phi a \}
\end{gather*}

The \textit{vertical composite}
\[
\begin{tikzcd}
	\A
		\arrow[r, "\Phi", ""{coordinate, name=A}, bend left=75, looseness=1.5, proarrow]
		\arrow[r, "\Psi", ""{coordinate, name=B}, proarrow]
		\arrow[r, "\Upsilon" swap, ""{coordinate, name=C}, proarrow, bend right=75, looseness=1.5]
		\arrow[from=A,to=B, phantom, "\Downarrow", "\scriptstyle\tau"{anchor=west, inner sep=0.8ex}, pos=0.4]
		\arrow[from=B,to=C, phantom, "\Downarrow", "\scriptstyle\sigma"{anchor=west, inner sep=0.8ex}]&
	\B
\end{tikzcd}
\]
in $\VCat^{\sharp}$ has component at $a \in \obj(\A)$ given by the composite morphism
\begin{gather*}
	\begin{tikzcd}[
		split arrows style,
		ampersand replacement=\&,
		column sep={between origins,12em},
	]
	I
		\xrightarrow{\,\tau_a\;}
	\sum_{x \in X} \A(a, x)
		\xrightarrow{\,\cong\;}
	\sum_{x \in X} \A(a, x) \otimes I
		\arrow[dr, phantom, ""{coordinate, name=Z}]
		\arrow[
			dr,
			rounded corners,
			"\,\sum \id \otimes \sigma_{x}"{anchor=south west, at start},
			to path={
				-- ([xshift=10ex]\tikztostart.east) \tikztonodes
				|- (Z)
				-| ([xshift=-2ex]\tikztotarget.west)
				-- (\tikztotarget)
			}
		]
		\&
	\\
		\&
	\sum_{x \in X} \A(a, x) \otimes 
		\big(\displaystyle\sum_{y \in Y} \A(x, y) \big)
		\xrightarrow{\,\cong\;}
	\sum_{x \in X} \A(a, x) \otimes \A(x, y)
		\xrightarrow{\,\sum[\mu_{a, x, y}]\;}
	\sum_{y \in Y} \A(a, y).
	\end{tikzcd}
	\\[1em]
	X = \Psi^{-1}\{\Phi a\}
	\qquad\qquad
	Y = \Upsilon^{-1}\{\Phi a\}
\end{gather*}

%----------------------------------------------------------------------%
% Section 3: The double catgeories                                     %
%----------------------------------------------------------------------%
\section{The double category of enriched categories, functors, and cofunctors}
\label{sec:double-category}

In this section, we introduce the double category $\VCOF$ 
whose objects are categories, whose horizontal morphisms are functors, and whose vertical morphisms are cofunctors, all enriched over a fixed distributive monoidal category $\V$. We reuse the name \textit{enriched natural transformation} for its cells because this notion subsumes both previously-defined notions of enriched natural transformation.

\begin{definition}
\label{definition:transformation-cell}
A \emph{enriched natural transformation}
\[
\begin{tikzcd}
\A
\arrow[d, proarrow, "\Phi"']
\arrow[r, "F"]
\arrow[dr, phantom, "\Rightarrow"{sloped}, "\scriptstyle\tau"{xshift=0.75ex,yshift=1.25ex}]
& \C
\arrow[d, proarrow, "\Psi"]
\\
\B
\arrow[r, "G"']
& \D
\end{tikzcd}
\]
of enriched functors and cofunctors consists of a family of morphisms 
in $\V$
\[
	\tau_{a} \colon I \longrightarrow 
	\displaystyle\sum_{x \in X} \C(Fa, x)\qquad\qquad X=\Psi^{-1}\{G\Phi a\}
\]
called the \textit{components} of $\tau$, such that the following diagram, encoding naturality, commutes.
\begin{gather*}
\begin{tikzcd}[column sep={3.3em}, cramped, ampersand replacement = \&, font=\small]
I \otimes \B(\Phi a, b)
\arrow[d, "\tau_{a} \otimes G_{\Phi a, b}"']
\& \B(\Phi a, b)
\arrow[l, "\cong"']
\arrow[r, "\cong"]
\&[0.6em] \B(\Phi a, b) \otimes I
\arrow[d, "\Phi_{a, b} \otimes \id"]
\\
\Big(\displaystyle\sum_{x \in X} \C(Fa, x) \Big) \otimes \D(G\Phi a, Gb)
\arrow[d, "\cong"']
\& \& \Big( \displaystyle\sum_{y \in Y} \A(a, y) \Big) \otimes I
\arrow[d, "\cong"]
\\
\displaystyle\sum_{x \in X} \C(Fa, x) \otimes \D(\Psi x, Gb)
\arrow[d, "\sum \id \otimes \Psi_{x, Gb}"']
\& \& \displaystyle\sum_{y \in Y} \A(a, y) \otimes I
\arrow[d, "\sum F_{a, y} \otimes \tau_{y}"]
\\
\displaystyle\sum_{x \in X} \C(Fa, x) \otimes 
\Big( \displaystyle\sum_{z \in Z} \C(x, z) \Big)
\arrow[d, "\cong"']
\& \& \displaystyle\sum_{y \in Y} \C(Fa, Fy) \otimes 
\Big( \displaystyle\sum_{z \in Z} \C(Fy, z) \Big)
\arrow[d, "\cong"]
\\
\displaystyle\sum_{z \in Z} 
\displaystyle\sum_{x \in X} \C(Fa, x) \otimes \C(x, z)
\arrow[r, "{\sum [\mu_{Fa, x, z}]}"', pos=0.44]
\& 
\displaystyle\sum_{z \in Z} \C(Fa, z)
\& \displaystyle\sum_{z \in Z} 
\displaystyle\sum_{y \in Y} \C(Fa, Fy) \otimes \C(Fy, z)
\arrow[l, "{\sum [\mu_{Fa, Fy, z}]}", pos=0.45]
\end{tikzcd}
\\[1em]
	X = \Psi^{-1}\{G\Phi a\} \qquad \qquad 
	Y = \Phi^{-1}\{b\} \qquad \qquad 
	Z = \Psi^{-1}\{Gb\}
\end{gather*}
\end{definition}

The \textit{horizontal composite}
\[
	\begin{tikzcd}
	\A
	\arrow[r, "F"]
	\arrow[d, proarrow, "\Phi"']
	\arrow[dr, phantom, "\Rightarrow"{sloped}, "\scriptstyle\tau"{xshift=0.75ex,yshift=1.25ex}]
	& \C
	\arrow[r, "H"]
	\arrow[d, proarrow, "\Psi"]
	\arrow[dr, phantom, "\Rightarrow"{sloped}, "\scriptstyle\sigma"{xshift=0.75ex,yshift=1.25ex}]
	& \mathcal{E}
	\arrow[d, proarrow, "\Upsilon"]
	\\
	\B
	\arrow[r, "G"']
	& \D
	\arrow[r, "K"']
	& \mathcal{F}
	\end{tikzcd}
\]
in $\VCOF$ has component at $a \in \obj(\A)$ given by the composite morphism
\begin{gather*}
	\begin{tikzcd}[
		split arrows style,
		ampersand replacement=\&,
		column sep={between origins,0em},
		row sep={0em},
		cramped
	]
	I
		\xrightarrow{\,\tau_{a}\;}
	\sum_{x \in X} \C(Fa, x)
		\xrightarrow{\,\sum H_{Fa, x}\;}
	\sum_{x \in X} \mathcal{E}(HFa, Hx)
		\xrightarrow{\,\cong\;}
	\sum_{x \in X} \mathcal{E}(HFa, Hx) \otimes I
		\arrow[dr, phantom, ""{coordinate, name=Z}]
		\arrow[
			dr,
			"\,\sum \id \otimes \sigma_{x}"{at start, anchor={south west}},
			rounded corners,
			to path={
				-- ([xshift=10ex]\tikztostart.east)\tikztonodes
				|- ([xshift=1ex]\tikztotarget.east)
				-- (\tikztotarget)
			}]
		\&[3em]
		\&[17em]
		\\
		\&
	\sum_{y \in Y}\sum_{x \in X} \mathcal{E}(HFa, Hx) \otimes \mathcal{E}(Hx, y)
		\xleftarrow{\;\cong\,}
	\sum_{x \in X} \mathcal{E}(HFa, Hx)  \otimes
		\big(\sum_{y \in Y} \mathcal{E}(Hx, y)\big)
		\arrow[dr, phantom, ""{coordinate, name=Y}]
		\arrow[
			dr,
			"{\sum [\mu_{HFa, Hx, y}]}\;"{swap, at end, anchor={north east}},
			rounded corners,
			to path={
				-- ([xshift=-2ex]\tikztostart.west) 
				|- ([xshift=-14ex]\tikztotarget.west)
				-- (\tikztotarget) \tikztonodes
			}]
		\&
		\\
		\&
		\&
	\sum_{y \in Y} \mathcal{E}(HFa, y).
	\end{tikzcd}
	\\[1em]
	X = \Psi^{-1}\{G\Phi a\}
	\qquad\qquad
	Y = \Upsilon^{-1}\{KG\Phi a\}
\end{gather*}

The \textit{vertical composite}
\[
\begin{tikzcd}
\A
\arrow[r, "F"]
\arrow[d, proarrow, "\Phi"']
\arrow[dr, phantom, "\Rightarrow"{sloped}, "\scriptstyle\tau"{xshift=0.75ex,yshift=1.25ex}]
& \C
\arrow[d, proarrow, "\Psi"]
\\
\B
\arrow[r, "G"]
\arrow[d, proarrow, "\Upsilon"']
\arrow[dr, phantom, "\Rightarrow"{sloped}, "\scriptstyle\sigma"{xshift=0.75ex,yshift=1.25ex}]
& \D
\arrow[d, proarrow, "\Omega"]
\\
\mathcal{E}
\arrow[r, "H"']
& \mathcal{F}
\end{tikzcd}
\]
in $\VCOF$ has component at $a \in \obj(\A)$  
given by the composite morphism
\begin{gather*}
	\begin{tikzcd}[
		split arrows style,
		ampersand replacement=\&,
		column sep={between origins,0em},
		row sep={0em},
		cramped
	]
	I
		\xrightarrow{\,\tau_{a}\;}
	\sum_{x \in X} \C(Fa, x)
		\xrightarrow{\,\cong\;}
	\sum_{x \in X} \C(Fa,x) \otimes I
		\xrightarrow{\,\sigma_{\Phi a}\;}
	\sum_{x \in X} \C(Fa, x) \otimes 
		\big(\sum_{y \in Y} \D(G\Phi a, y) \big)
		\arrow[dr, phantom, ""{coordinate, name=Z}]
		\arrow[
			dr,
			"\,\cong"{at start, anchor={south west}},
			rounded corners,
			to path={
				-- ([xshift=8.5ex]\tikztostart.east)\tikztonodes
				|- ([xshift=1ex]\tikztotarget.east)
				-- (\tikztotarget)
			}]
		\&[3em]
		\&[7em]
	\\
		\&
	\sum_{y \in Y}\sum_{x \in X} \C(Fa, x) 
		\otimes \big(\sum_{z \in Z_{y}} \C(x, z) \big)
		\xleftarrow{\;\sum \sum \id \otimes \Psi_{x, y}\,}
	\sum_{y \in Y}\sum_{x \in X}
		\C(Fa, x) \otimes \D(\Psi x, y)
		\arrow[dr, phantom, ""{coordinate, name=Y}]
		\arrow[
			dr,
			"\cong\;"{at end, anchor={south east}},
			rounded corners,
			to path={
				-- ([xshift=-2ex]\tikztostart.west) 
				|- ([xshift=-9ex]\tikztotarget.west)
				-- (\tikztotarget) \tikztonodes
			}]
		\&
	\\
		\&
		\&
	\sum_{w \in W} \displaystyle\sum_{x \in X}
		\C(Fa, x) \otimes \C(x, w)
		\xrightarrow{\,\sum [\mu_{Fa, x, w}]\;}
	\sum_{w \in W} \C(Fa, w).
	\end{tikzcd}
	\\[1em]
	X = \Psi^{-1}\{G\Phi a\} \qquad
	Y = \Omega^{-1}\{ H\Upsilon\Phi a\} \qquad
	Z_{y} = \Psi^{-1}\{y\} \qquad 
	W = (\Omega\Psi)^{-1}\{H\Upsilon\Phi a \}
\end{gather*}

Recall that every double category has a horizontal and a vertical $2$-category \cite{GP99}. 

\begin{proposition}
The horizontal and vertical $2$-categories of $\VCOF$ are
$\VCat$ and $\VCat^{\sharp}$.
\end{proposition}

%----------------------------------------------------------------------%
% Section 5: Enriched lenses                                           %
%----------------------------------------------------------------------%
\section{Enriched lenses as compatible squares}
\label{sec:enriched-lenses}

In this section, we introduce the notion of enriched lens as an 
enriched functor and enriched cofunctor that satisfy a certain 
compatibility condition. 
We also introduce the flat double subcategory $\VCOF_{\flt}$ of 
$\VCOF$, whose cells are called \emph{compatible squares}, 
and show that enriched lenses arise naturally as certain cells in
this double category. 

\begin{definition}
\label{definition:enriched-lens} 
An \emph{enriched lens} $(F, \Phi) \colon \A \rightleftharpoons \B$ 
consists of an enriched functor $F \colon \A \rightarrow \B$ and an 
enriched cofunctor $\Phi \colon \A \pto \B$, such that 
$F a = \Phi a$ for all $a \in \obj(\A)$ and the following diagram commutes.
\begin{equation}
\label{equation:enriched-lens}
\begin{tikzcd}[column sep = tiny]
\B(Fa, b)
\arrow[rd, "\Phi_{a, b}" swap]
\arrow[rr, "\id"]
& & \B(Fa, b)
\\
& \displaystyle\sum_{x \in X} \A(a, x) 
\arrow[ru, "{[F_{a, x}]}" swap]
\end{tikzcd}
\qquad \qquad X = F^{-1}\{b\}
\end{equation}
\end{definition}

There is a category $\VLens$ whose objects are enriched categories and whose morphisms are enriched lenses, 
where composition of enriched lenses is defined by the composition of 
the underlying enriched functors and enriched cofunctors. 

\begin{remark}
The category $\VLens$ also has a $2$-category structure. 
A $2$-cell between a pair of enriched lenses consists of 
a enriched natural transformation between their underlying functors 
and an enriched natural transformation between their underlying cofunctors
that are compatible with each other (in a similar sense to 
Definition~\ref{definition:compatible-square} below).
\end{remark}

\begin{definition}
\label{definition:compatible-square}
A square of functors and cofunctors
\begin{equation}
	\label{equation:compatible-square}
	\begin{tikzcd}
	\A
	\arrow[d, proarrow, "\Phi"']
	\arrow[r, "F"]
	& \C
	\arrow[d, proarrow, "\Psi"]
	\\
	\B
	\arrow[r, "G"']
	& \D
	\end{tikzcd}
\end{equation}
is a \emph{compatible square} if $G\Phi a = \Psi F a$, and the following diagram commutes.
\begin{equation*}
\begin{tikzcd}[column sep = large]
\B(\Phi a, b)
\arrow[d, "\Phi_{a, b}"']
\arrow[r, "G_{\Phi a, b}"]
& \D(\Psi F a, Gb)
\arrow[d, "\Psi_{Fa, Gb}"]
\\
\displaystyle\sum_{x \in X} \A(a, x)
\arrow[r, "\sum F_{a, x}"']
& \displaystyle\sum_{y \in Y} \C(Fa, y)
\end{tikzcd}
\qquad \qquad 
\begin{aligned}
X &= \Phi^{-1}\{b\} \\
Y &= \Psi^{-1}\{Gb\}
\end{aligned}
\end{equation*}
\end{definition}

\begin{remark}
	We should think of compatibility as a commutativity-like notion for diagrams involving both functors and cofunctors of this particular shape.
\end{remark}

There is a double category $\VCOF_{\flt}$ whose objects are 
enriched categories, whose horizontal morphisms are enriched functors, 
whose vertical morphisms are enriched cofunctors, and whose 
cells are compatible squares. 
This is a double subcategory of $\VCOF$ via the double functor 
that sends a compatible square \eqref{equation:compatible-square} 
to the enriched natural transformation with the same boundary 
whose component at $a \in \obj(\A)$ is given by
the composite below. 
\begin{equation*}
\begin{tikzcd}
I 
\arrow[r, "\eta_{Fa}"]
& \C(Fa, Fa) 
\arrow[r, "\iota_{Fa}"]
& \displaystyle\sum_{x \in X} \C(Fa, x) 
\end{tikzcd}
\qquad \qquad
X = \Psi^{-1}\{ G\Phi a \}
\end{equation*}

\begin{proposition}
\label{proposition:lens-as-compatible-square}
An enriched lens $(F, \Phi) \colon \A \rightleftharpoons \B$ is 
precisely a compatible square
\begin{equation*}
\begin{tikzcd}[cramped]
\A
\arrow[d, proarrow, "\Phi"']
\arrow[r, "F"]
& \B
\arrow[d, proarrow, "1_{\B}"]
\\
\B
\arrow[r, "1_{\B}"']
& \B
\end{tikzcd}
.
\end{equation*}
\end{proposition}

%----------------------------------------------------------------------%
% Section 6: Examples                                                  %
%----------------------------------------------------------------------%
\section{Examples of enriched cofunctors and lenses}
\label{sec:examples}

In this section, we describe the enriched cofunctors and lenses for several interesting bases of enrichment.

%----------------------------------------------------------------------%
% Subsection                                                           %
%----------------------------------------------------------------------%
\subsection{Locally small categories} 
Consider the category $\Set$ equipped with the cartesian monoidal 
product.
This is a distributive monoidal category, with coproducts given by 
disjoint union. 
Categories, cofunctors, and lenses enriched in $\Set$ are precisely 
\emph{locally small categories}, 
\emph{cofunctors} \cite{Agu97},
and \emph{delta lenses} \cite{DXC11, JR13}, respectively. 

A cofunctor $\Phi \colon \A \pto \B$ consists of a 
function $\Phi \colon \obj(\A) \rightarrow \obj(\B)$ and a 
family
\begin{equation}
\label{equation:lifting-operation}
\Phi_{a, b} \colon \B(\Phi a , b) 
\longrightarrow \!\!\!\! \sum_{x \in \Phi^{-1}\{b\}} \!\!\!\! \A(a, x) 
\end{equation}
of lifting maps such that $\Phi_{a, \Phi a}(1_{\Phi a}) = 1_{a}$ and 
$\Phi_{a, b'}(v \circ u) = \Phi_{a', b'}(v) \circ \Phi_{a, b}(u)$.
The equation
$\Phi(\cod\Phi_{a, b}(u)) = \cod u = b$ holds by construction. 
Examples of cofunctors $\A \pto \B$ include 
bijective-on-objects functors $\B \rightarrow \A$ and discrete 
opfibrations $\A \rightarrow \B$.

A delta lens $(F, \Phi) \colon \A \rightleftharpoons \B$ consists 
of a functor $F \colon \A \rightarrow \B$ and a cofunctor 
$\Phi \colon \A \pto \B$ such that $F a = \Phi a$, 
and $F_{a,a'}\Phi_{a, b}(u) = u$ where 
$a' = \cod \Phi_{a, b}(u)$. 
Split opfibrations are precisely delta lenses $(F, \Phi)$
where each morphism
$\Phi_{a, b}(u)$ is an opcartesian lift along $F$. 

A natural transformation
\[
\begin{tikzcd}
\A
\arrow[d, cofunctor, "\Phi"']
\arrow[r, "F", functor]
\arrow[dr, phantom, "\Rightarrow"{sloped}, "\scriptstyle\tau"{xshift=0.75ex,yshift=1.25ex}]
& \C
\arrow[d, cofunctor, "\Psi"]
\\
\B
\arrow[r, "G"', functor]
&
\D
\end{tikzcd}
\]
of functors and cofunctors is a family 
\[\tau_a \colon Fa \to \cod \tau_a \qquad\qquad \Psi \cod \tau_a = G\Phi a\]
of morphisms of $\C$ such that, for each $u \colon b \to b'$ in $\B$ and each $a \in \Phi^{-1}\{b\}$, the square
\begin{gather*}
\begin{tikzcd}[ampersand replacement=\&]
Fa
	\arrow[d, "F \Phi_{a,b'}(u)" swap]
	\arrow[r, "\tau_a"]
	\&
c
	\arrow[d, "\Psi_{b,Gb'}G(u)"]
\\
Fa'
	\arrow[r, "\tau_{a'}" swap]
	\&
c'
\end{tikzcd}
\\[1em]
a' = \cod \Phi_{a,b'}(u)
	\qquad\qquad
c  = \cod \tau_a
	\qquad\qquad
c' = \cod \tau_{a'} = \cod \Psi_{b,Gb'}G(u)
\end{gather*}
in $\C$, called the \textit{naturality square} of $\tau$ at $u$, always commutes.

%----------------------------------------------------------------------%
% Subsection                                                           %
%----------------------------------------------------------------------%
\subsection{Linear categories} 
\label{sec:linear-categories}
Consider the category $\Vec_{K}$ of vector spaces over the field $K$,
equipped with the usual tensor product of vector spaces. 
This is a distributive monoidal category, with 
coproducts given by the direct sum of vector spaces. 

A category enriched in $\Vec_{K}$ is called a 
\emph{linear category} or \emph{algebroid} 
(an algebroid with one object is an algebra over a field).  
Aguiar showed that 
every category enriched in $\Vec_{K}$ may be viewed as a 
category internal to $\Vec_{K}$ 
(using a more general definition of category internal to a monoidal 
category with nice equalisers)~\cite[Section~9.1]{Agu97};
we now observe that the internal cofunctors between a pair of 
enriched categories, via this identification, 
are exactly the enriched cofunctors between these enriched categories.

A cofunctor enriched in $\Vec_{K}$ will be called a \textit{linear cofunctor}.
Recall that the direct sum of a family of vector spaces is the subspace of their direct product formed by those tuples that are non-zero in only finitely many coordinates. Hence, a linear cofunctor $\Phi \colon \A \pto \B$ chooses, for each morphism $u \colon \Phi a \rightarrow b$ in $\B$, not just one lift, but a family
\[\big\{ \Phi_{a,b}(u)_x \colon a \to x \big\}_{x \in X_{a,b,u}}\]
of lifts, indexed by a finite subset $X_{a,b,u}$ of $\Phi^{-1} \{b\}$. Of course, these choices are subject to the requirement that the lifting maps $\Phi_{a,b}$ be linear. As finite direct sums are also categorical products, when the fibre $\Phi^{-1}\{b\}$ is finite, the data of the lifting map
\[\Phi_{a,b} \colon \B(\Phi a, b) \longrightarrow \!\!\!\!\bigoplus_{x \in \Phi^{-1}\{b\}} \!\!\!\! \A(a, x)\]
is equivalent to the data of a family of linear maps
\[\big\{\Phi_{a,b,x} \colon \B(\Phi a, b) \longrightarrow \A(a, x) \big\}_{x \in \Phi^{-1}\{b\}}.\]

Recall that linear categories and functors may be viewed as (ordinary) categories and functors via base change along the forgetful functor $\Vec_{K} \rightarrow \Set$; this is possible because the forgetful functor is lax monoidal. As the forgetful functor does not preserve coproducts, the same is not true of linear cofunctors. In contrast, the free functor $\Set \rightarrow \Vec_{K}$ 
is both strong monoidal \emph{and} preserves coproducts, so base change along the free functor works not only for categories and functors, but also for cofunctors.

%----------------------------------------------------------------------%
% Subsection                                                           %
%----------------------------------------------------------------------%
\subsection{Lawvere metric spaces} 
\label{sec:Lawvere-metric-spaces}

Consider the thin category $([0, \infty], \geq)$ of extended 
non-negative real numbers, equipped with addition
$ + \colon [0, \infty] \times [0, \infty] \rightarrow [0, \infty]$
as its monoidal product 
and the terminal object $0$ as its monoidal unit. 
This is a distributive monoidal category, 
with coproducts given by infima (greatest lower bounds). 

A category $\A$ enriched in $([0, \infty], \geq)$ is a 
\emph{Lawvere metric space} \cite{Law73}. 
It consists of a set of objects $\obj(\A)$ together with a 
real number $\dist(a, a') \in [0, \infty]$ for each pair of objects,
which we think of as the \emph{distance} from $a$ to $a'$. 
The existence of the identity elements and composition maps of $\A$
is equivalent to requiring that
\[
	\dist(a, a) = 0
	\qquad\text{and}\qquad
	\dist(a, a') + \dist(a', a'') \geq \dist(a, a'').
\]

An functor $F \colon \A \rightarrow \B$ enriched in $([0, \infty], \geq)$ is a function 
$F \colon \obj(\A) \rightarrow \obj(\B)$ with
\[\dist(a, a') \geq \dist(Fa, Fa').\]
Such functions are called \textit{non-expansive} or \textit{1-Lipschitz}. 

A cofunctor $\Phi \colon \A \pto \B$ enriched in $([0, \infty], \geq)$ is a function 
$\Phi \colon \obj(\A) \rightarrow \obj(\B)$ 
with
\begin{equation}
	\label{Lawvere cofunctor}
	\dist(\Phi a, b) \geq \inf_{x \in \Phi^{-1}\{b\}} \dist(a, x).
\end{equation}
This notion is closely related to the notion of \textit{submetry}.\footnote{Many thanks to Callum Reader for pointing us towards this notion.}

A \textit{submetry} $\Phi \colon \A \to \B$ between Lawvere metric spaces is a function $\Phi \colon \obj(\A) \to \obj(\B)$ such that for each $a \in \obj(\A)$ and each $b \in \obj(\B)$, there is an $x \in \Phi^{-1}\{b\}$ such that $d(\Phi a, b) \geq d(a, x)$. Although every submetry is an enriched cofunctor, not every enriched cofunctor is a submetry. For example, consider the Lawvere metric spaces~$\A$ and $\B$ where
\[\obj(\A) = \{0\} \cup (1, 2],\qquad \obj(\B) = \{0\} \cup \{1\},\quad\text{and}\quad
	d(x, y) = \begin{cases}
		y - x &\text{if } x \leq y,\\
		\infty &\text{otherwise.}
	\end{cases}
\]
The enriched cofunctor $\Phi \colon \A \pto \B$ with object map
\[
	\Phi(a) = \begin{cases}
		0 &\text{if } a = 0,\\
		1 &\text{otherwise}
	\end{cases}
\]
is not a submetry because $d(\Phi(0), 1) = 1$ whilst $d(0, x) = x > 1$ for each $x \in \Phi^{-1} \{1\}$.

An enriched lens $F \colon \A \rightleftharpoons \B$ is a 
function $F \colon \obj(\A) \rightarrow \obj(\B)$ such that
\[
	\dist(Fa, b) = \inf_{x \in F^{-1}\{b\}} \dist(a, x).
\]
Even an enriched lens is not necessarily a submetry, indeed the enriched cofunctor $\Phi$ defined in the example above is actually an enriched lens.

\begin{remark}
If a cofunctor $\Phi \colon \A \to \B$ enriched in $([0,\infty], \geq)$ has \textit{finite} fibres, then the infimum in \eqref{Lawvere cofunctor} is always attained, and so $\Phi$ is actually a submetry. 
Finding conditions on a cofunctor that are both sufficient and \textit{necessary} to make it a submetry is an open problem.
\end{remark}

%----------------------------------------------------------------------%
% Subsection                                                           %
%----------------------------------------------------------------------%
\subsection{Enrichment in the free coproduct completion of a monoidal category}
\label{sec:Fam-construction}

One aspect where enriched cofunctors differ from enriched functors 
is the requirement for $\V$ to be a 
\emph{distributive monoidal category}. 
Whilst it is not uncommon for many monoidal categories to be distributive, 
one might ask if it is possible to describe enriched cofunctors starting 
with an arbitrary monoidal category $(\W, \otimes, I)$. 

The \emph{free coproduct completion} $\Fam(\C)$ 
of a category $\C$ is the comma category below. 
\begin{equation*}
\begin{tikzcd}[column sep=small]
\Fam(\C) 
\arrow[rd, phantom, "\Rightarrow"]
\arrow[r]
\arrow[d]
& \ast 
\arrow[d, "\C"]
\\
\Set
\arrow[r, hook, "\disc"']
& \Cat
\end{tikzcd}
\end{equation*}
In detail, the objects of $\Fam(\C)$ are \emph{families} of objects in $\C$, 
that is, functions $J \rightarrow \obj(\C)$ or, equivalently,
functors $\disc(J) \rightarrow \C$ where $\disc(J)$ is the 
discrete category on the set $J$. 
A morphism in $\Fam(\C)$ from $j \colon \disc(J) \rightarrow \C$ to 
$k \colon \disc(K) \rightarrow \C$ consists of a function 
$f \colon J \rightarrow K$ together 
with a natural transformation as depicted below.
\begin{equation*}
\begin{tikzcd}[column sep = tiny]
\disc(J)
\arrow[rr, "\disc(f)"]
\arrow[rd, "j"']
& 
\phantom{X}
\arrow[d, phantom, "\Rightarrow", pos=0.3]
& \disc(K)
\arrow[ld, "k"]
\\
& \C & 
\end{tikzcd}
\end{equation*}
The coproduct of objects $j \colon \disc(J) \rightarrow \C$ and 
$k \colon \disc(K) \rightarrow \C$ is given by
\begin{equation*}
\begin{tikzcd}
\disc(J \sqcup K) \cong \disc(J) + \disc(K) 
\arrow[r, "{[j, k]}"]
& \C
\end{tikzcd}
\end{equation*}
and the initial object is given by the unique functor 
from the empty category to $\C$. 

For a monoidal category $(\W, \otimes, I)$, the category
$\Fam(\W)$ is a distributive monoidal. 
The monoidal product of objects $j \colon \disc(J) \rightarrow \C$ and 
$k \colon \disc(K) \rightarrow \C$ is given by 
\begin{equation*}
\begin{tikzcd}
\disc(J \times K) \cong \disc(J) \times \disc(K) 
\arrow[r, "j \times k"]
& \W \times \W
\arrow[r, "\otimes"]
& \W
\end{tikzcd}
\end{equation*}
with unit given by the functor from the terminal category to $\W$ 
that chooses the object~$I$. 
The distributivity of $\Fam(\W)$ follows from the distributivity of 
$\Set$ equipped with the cartesian monoidal product.  

For simplicity, we consider the case where $\W$ is a thin category. 
We will write $x \geq y$ if and only if a morphism $x \rightarrow y$ exists.

A category enriched in $\Fam(\W)$ consists of a locally small category
$\A$ together with an object $|u|$ in $\W$ for each morphism 
$u \colon a \rightarrow a'$ in $\A$, 
such that
\begin{equation}
\label{equation:Fam-category}
	I \geq |1_{a}| \qquad \text{and} \qquad
	|u| \otimes |v| \geq |v \circ u|.
\end{equation}
A functor enriched in $\Fam(\W)$ is a functor 
$F \colon \A \rightarrow \B$ such that
\begin{equation}
\label{equation:Fam-functor}
	|w| \geq |F_{a, a'}(w)|.
\end{equation}
A cofunctor enriched in $\Fam(\W)$ is a cofunctor 
$\Phi \colon \A \pto \B$ such that 
\begin{equation}
\label{equation:Fam-cofunctor}
	|u| \geq |\Phi_{a, b}(u)|.
\end{equation}
A lens enriched in $\Fam(\W)$ is delta lens 
$(F, \Phi) \colon \A \rightleftharpoons \B$ such that 
both \eqref{equation:Fam-functor} and \eqref{equation:Fam-cofunctor} hold;
by \eqref{equation:enriched-lens}, it follows that $|u| = |\Phi_{a, b}(u)|$
since $\W$ is a thin category.

%----------------------------------------------------------------------%
% Subsection                                                           %
%----------------------------------------------------------------------%
\subsection{Weighted categories}
\label{sec:weighted-lenses}

Let $\wSet$ be the free coproduct completion $\Fam([0, \infty], \geq)$
of the poset $([0, \infty], \geq)$ of extended non-negative real numbers 
(see Sections~\ref{sec:Lawvere-metric-spaces} and \ref{sec:Fam-construction}). 
The objects of this category are called \emph{weighted sets} \cite{Gra07}. 
Perrone introduced the notion of \emph{weighted lenses}~\cite{Per21}, 
which we now recognise as lenses enriched in $\wSet$.

A category enriched in $\wSet$ consists a locally small category $\A$ 
together with a \emph{weight} $|u| \in [0, \infty]$ for each morphism
$u \colon a \rightarrow a'$ in $\A$, such that the inequalities
\eqref{equation:Fam-category} hold. 
As the monoidal unit $0$ is also a terminal object, 
the first inequality in \eqref{equation:Fam-category} may be strengthened
to the equality $|1_{a}| = 0$.

Weighted functors are functors such that the
condition \eqref{equation:Fam-functor} holds. 
Similarly, weighted cofunctors are cofunctors 
such that the condition \eqref{equation:Fam-cofunctor} holds.
Weighted lenses are delta lenses such that the conditions 
\eqref{equation:Fam-functor} and \eqref{equation:Fam-cofunctor} hold.
Since $[0, \infty]$ is a poset, we may deduce the stronger condition
$|\Phi_{a, b}(u)| = |u|$ on the weight of a lift. 

%----------------------------------------------------------------------%
% Section                                                           %
%----------------------------------------------------------------------%
\newpage
\section*{Acknowledgements}
This collaboration between the two authors began when they
were research students at Macquarie University. 
The resulting work was first presented at the Australian Category Seminar 
in October 2021, 
and later at the Applied Category Theory conference in July~2022. 
We would like to thank everyone who provided useful feedback and ideas
on this research, with a special thanks to Paolo Perrone and 
Michael Johnson. 

%-----------------------------------------------------------------------
% References 
%-----------------------------------------------------------------------
\raggedright
\printbibliography

\end{document}